\documentclass[12pt,a4paper]{article}
\usepackage{amssymb,amsmath,amsthm}
\newtheorem{theorem}{\indent {\sc Theorem}}
\begin{document}

{\large

\begin{center}
\textbf{\Large{Hua Loo Keng's Problem Involving Primes\\ of a
Special Type}}
\end{center}

\bigskip

\begin{center}
\textbf{\large{Sergey A. Gritsenko and Natalya N. Motkina}}
\end{center}
\bigskip
\begin{abstract}
Let $\eta$ be a quadratic irrationality. The variant of Hua Loo
Keng's problem  involving primes such that $a<\{\eta p^2\}<b$,
where $a$ and $b$ are arbitrary real numbers of the interval
$(0,1)$, solved in this paper.
\end{abstract}
\bigskip

\textbf{1. Introduction} Suppose that $k\ge 2$ and $n\ge 1$ are
positive integers. Consider the equation
\begin{equation}\label{e1}
p_1^n+p_2^n+\cdots +p_k^n=N
\end{equation}
in primes $p_1,\ p_2,\ldots ,p_k$.

Several classical problems in number theory can be reduced to the
question on the number of solutions of the equation  (\ref{e1}).
For example, if $n=1$ and $k=2$ or $3$ then (\ref{e1}) is the
equation of Goldbach; if $n\ge 3$ then (\ref{e1}) is the equation
of Waring--Goldbach.

Let $\mathcal{P}$ be a subset of the set of primes,
$\mu=\lim\limits_{n\to\infty}\pi^{-1}(N)\sum\limits_{\substack{p\le
N\\p\in \mathcal{P}}}1$ be the "density" of $\mathcal{P}$,
$0<\mu<1$.

Let $J_{k,n}(N)$ be the number of solutions of (\ref{e1}) in
primes of the set $\mathcal{P}$, and $I_{k,n}(N)$ be the number of
solutions of (\ref{e1}) in arbitrary primes.

Sometimes the equality
\begin{equation}\label{e2}
J_{k,n}(N)\sim\mu^k I_{k,n}(N)
\end{equation}
holds.

For example, if $\mathcal{P}=\{p\ |\ a<\{\sqrt{p}\}< b \}$,
$0<a<b<1$ then provided $n=1$, $k=3$, and also provided $n\ge 3$,
$k\gg n^2\log n$ the equality (\ref{e2}) holds (s. \cite {G}).

In present paper we solve the additive problem for which the
pro\-per\-ty (\ref{e2}) does not hold.

Let $\eta$ be a quadratic irrationality, $0<a<b<1$,
$\mathcal{P}=\{p\ |\ a<\{\eta p^2\}<b\}$, $n=1$. Then $\mu=b-a$
(s. \cite {V}).

The main result of the present paper is the following theorem
which shows that the equality
\begin{equation}\label{e3}
J_{5,2}(N)\sim (b-a)^5 I_{5,2}(N)
\end{equation}
does not hold.

Note that (s. \cite {Vv})
$$ I_{5,2}(N)\asymp \frac{N^{3/2}}{\log^5 N}
$$
provided $N\equiv 5 (\mod 24)$.

\begin{theorem}
Suppose
$$
\sigma(N,a,b)=\sum_{|m|<\infty} e^{2\pi i m(\eta N-2.5(a+b))}
\frac{\sin^5 \pi m (b-a)}{\pi ^5 m^5}.
$$
Then the following formula holds
$$
J_{5,2}(N)=I_{5,2}(N)\sigma(N,a,b)+O(N^{3/2-0,00002}).
$$
\end{theorem}

{\textbf{Outline of proof.}  Define the function
$$\psi_0 (x)= \bigl\{
\begin{array}{ll} 1,\quad \emph{if}\quad & a<x<b,\\ 0,
\quad \emph{if}\quad & 0\leq x\leq a \quad \emph{or} \quad b\leq x
\leq 1
\end{array} $$
and extend it by periodicity to the entire numerical axis. So that
$$ J_{5,2}(N)=\int_0^1 \left(\sum_{p\leq \sqrt
N}\psi_0 (\eta p^2) e^{2 \pi i x p^2}\right)^5 e^{-2 \pi i x N}dx.
$$
We record  Vinogradov's lemma  (see [1], Section III, Chapter 2,
Lemma 2). In the notation of the lemma, let $r=[\ln N]$, $\Delta =
N^{-0,01}$. Denote by $\psi_1$ a function $\psi$ from the lemma
with $\alpha=a+\Delta/2$, $\beta=b-\Delta/2$ and by $\psi_2$ with
$\alpha=a-\Delta/2$, $\beta=b+\Delta/2$. Then we have the
inequality
\begin{equation}\label{e4}
J_1(N)\leq J_{5,2}(N) \leq J_2(N),
\end{equation}
where
$$ J_k(N)=\int_0^1 \left(\sum_{p\leq \sqrt N}\psi_k (\eta p^2)
e^{2 \pi i x p^2}\right)^5 e^{-2 \pi i x
N}dx,\quad k=1,2.
$$
Next we derive an asymptotic formulas for $J_1(N)$ and $J_2(N)$.
The main terms of these formulas coincide. Use
the Fourier series for $\psi_k(x)$
$$
\psi_k (\eta p^2)=\sum_{|m|\leq r\Delta^{-1}}c_k(m)e^{2 \pi i m
\eta p^2}+O(N^{-\ln \pi}),
$$
thus
$$
J_k(N)=\sum_{|m_1|\leq r\Delta ^ {-1}} c_k(m_1) \sum_{|m_2|\leq
r\Delta^{-1}}c_k(m_2)\cdot$$
$$\cdot\sum_{|m_3|\leq r\Delta^{-1}}c_k(m_3)
\sum_{|m_4|\leq r\Delta^{-1}}c_k(m_4)\sum_{|m_5|\leq r\Delta^{-1}}
c_k(m_5)\cdot$$
$$\cdot\int_0^1 S(x+
m_1\eta)S(x+ m_2\eta)S(x+ m_3\eta)S(x+ m_4\eta)S(x+ m_5\eta)e^{-2
\pi i x N}dx$$ $$+O(N^{3/2-\ln \pi}\ln N),
$$
where
$$S(x)=\sum_{p\leq \sqrt{N}}e^{2\pi i x p^2}.
$$

Note that
$$
\sum_{|m|\leq r\Delta^{-1}}c^5_k(m)\int_0^1 S^5(x+m\eta)e^{-2\pi
ixN}dx=
$$
\begin{equation}\label{e5}
=I_{5,2}(N)(\sigma(N,a,b)+O(\Delta)).
\end{equation}
Consider the sets such that $(m_1,m_2,m_3,m_4,m_5)\neq(m,m,m,m,m)$
and the integrals
$$I(N,m_1,m_2,m_3,m_4,m_5)=$$ $$=\int_0^1 |S(x+ m_1\eta)||S(x+
m_2\eta)||S(x+ m_3\eta)||S(x+m_4\eta)||S(x+ m_5\eta)|dx.
$$
Without loss of generality it can be assumed that $m_1<m_2$. By
$t$ denote  $x+ m_1 \eta$. Since the integrand has a period $1$,
we can take $t$ to lie in the range $E=[-1/\tau;1-1/\tau)$, where
$\tau=N^{1-0,001}$. Suppose  $d,q\in \mathbb{Z}$ are such that
$$
t=\frac{d}{q}+\frac{\theta_1}{q\tau},\quad (d,q)=1,\quad 1\leq
q\leq \tau,\quad |\theta_1|<1.
$$
Divide the interval $E$ into two intervals $E_1$ and $E_2$ for
which $q\leq N^{0,001}$ and $N^{0,001}<q\leq \tau$ respectively.
On writing $m'_2=m_2-m_1$, $\ldots$, $m'_5=m_5-m_1$, we can
express the above integral as
$$
\int_{-1/\tau}^{1-1/\tau} |S(t)||S(t+m'_2\eta)||S(t+m'_3
\eta)||S(t+m'_4 \eta)||S(t+m'_5 \eta)|dt=$$
$$
=\int_{E_1}F(t)dt+\int_{E_2}F(t)dt,
$$
where
$$
F(t)=|S(t)||S(t+m'_2\eta)||S(t+m'_3 \eta)||S(t+m'_4
\eta)||S(t+m'_5 \eta)|.$$ It is well known that if $t\in E_2$ then
$$ |S(t)|\ll N^{0,5-0,00002}.$$
Consider $|S(t+m'_2\eta)|$ when $t\in E_1$. Since $\eta$ is the
quadratic irrationa\-li\-ty then we have
\begin{equation}\label{e6}
t+m'_2\eta=\frac{X}{Y}+\frac{\theta_2}{Y^2}, \quad (X,Y)=1,\quad
|\theta_2|<N^{0,012}\ln N
\end{equation} and
\begin{equation}\label{e7}
\sqrt{\tau} N^{-0,01015}\ln N <Y< \sqrt{\tau} N^{0.001}.
\end{equation}
Then from (\ref{e6}) and (\ref{e7}) it follows in the standard way
the estimate
$$ |S(t+m'_2\eta)|\ll N^{0,5-0,00002}.$$
Using the inequality between arithmetical and geometric means, we
finally obtain
\begin{equation}\label{e8}
I(N,m_1,m_2,m_3,m_4,m_5)\ll N^{3/2-0,00002}.
\end{equation}
Now the theorem  follows directly from formulas (\ref{e4}),
(\ref{e5}), (\ref{e8}).

\textit{Remark.} If $\mathcal{P}=\{p\ |\ a<\{\eta p\}<b\}$ then
the equality (\ref{e3}) holds.

\end{document}